\documentclass[12pt]{amsart}

\usepackage{enumerate}
 
\usepackage[active]{srcltx}
\usepackage{nicefrac}

\usepackage{amsthm,amssymb,setspace}
\usepackage[pagebackref,colorlinks,linkcolor=red,citecolor=blue,urlcolor=blue,hypertexnames=true]{hyperref}
\usepackage{amsrefs}
\usepackage[matrix, arrow]{xy}

\newcommand{\N}{\mathbb N}

\newcommand{\cC}{\mathcal C}
\newcommand{\cF}{\mathcal F}
\newcommand{\cS}{\mathcal S}
\newcommand{\cU}{\mathcal U}

\newcommand{\F}{\mathbb F}

\newcommand{\Z}{\mathbb Z}

\theoremstyle{plain}
\newtheorem{theorem}{Theorem}[section]
\newtheorem{corollary}[theorem]{Corollary}
\newtheorem{lemma}[theorem]{Lemma}
\newtheorem{proposition}[theorem]{Proposition}

\theoremstyle{definition}

\newtheorem{definition}[theorem]{Definition}

\theoremstyle{remark}

\newtheorem{remark}[theorem]{Remark}

\begin{document}

\onehalfspace

\title{About the metric approximation of Higman's group}

\author{Andreas Thom}
\address{Andreas Thom, Universit\"at Leipzig, Germany}
\email{thom@math.uni-leipzig.de}

\begin{abstract}
We prove that Higman's group does not embed into a metric ultraproduct of finite groups with a commutator-contractive invariant length function.
\end{abstract}

\maketitle

%\tableofcontents

\section{Introduction}

This article is about the metric approximation of group laws by manageable structures. A first occurance of this kind of question was in Alain Connes' seminal work \cite{MR0454659} where he noted that the group von Neumann algebra of the free group can be well-approximated by matrix-algebras, even though it is not itself hyperfinite. He conjectured, that this ought to be true for all tracial von Neumann algebras.
In recent years, mostly through the work of Eberhard Kirchberg \cite{MR1218321}, the Connes embedding problem for tracial von Neumann algebras has attracted again a lot of attention; see \cite{MR2072092} for various equivalent formulations, references and recent progress. Also, a $C^*$-algebraic analogue of the question has been studied in detail, see \cites{MR1437044, MR1796908, MR1835511}. From various points of view, these problems seem to be notoriously difficult. The particular case of group von Neumann algebras has been studied deeply and analogous questions in a more combinatorial setup have been asked as well. Misha Gromov's work \cite{MR1694588} on Gottschalk's Surjunctivity Conjecture \cite{MR0407821} has lead to the notion of \emph{soficity} for discrete group (see \cites{MR2178069, MR2089244}) and it remains to be an open problem whether or not all groups are sofic. The importance of this question is due to the fact that various conjectures about discrete groups could be proved for sofic groups \cite{MR2178069, MR2089244,MR2417890}.
A common theme which has been singled out is the metric approximation of group laws by certain classes of finite or compact groups with a specified invariant metric.

Before start to explain what we mean, we have to fix some conventions and definitions. We use the notation
$[g,h] = ghg^{-1}h^{-1}$ and $\bar g = g^{-1}$.

\begin{definition}
Let $G$ be a group. An invariant length function is a map
$\ell \colon G \to [0,1]$ such that $\ell(g)=0$ if and only if $g=e$, and
$$\ell(gh) \leq \ell(g) + \ell(h), \quad \ell(g^{-1}) = \ell(g), \quad \mbox{and} \quad
\ell(hg) = \ell(gh)$$
for all $g,h \in G$.
\end{definition}

Note that the third condition is equivalent to the requirement $\ell(g)=\ell(hg\bar h)$ for all $g,h \in G$. Every invariant length function induces a metric
$d(g,h) := \ell(g\bar h)$ on $G$ which satisfies 
$$d(kg,kh) = d(gk,hk) = d(g,h), \quad \forall g,h,k \in G.$$

\begin{remark}
It would have been natural to choose $[0,2]$ instead of $[0,1]$ from the operator-algebraic point of view, but we will stick to a convention which is more compatible with the point of view of model theory for metric structures, see \cite{MR2436146}. 
\end{remark}

The following condition on an invariant length function will be crucial in this note.
\begin{definition}
An invariant length function is said to be \emph{commutator-contractive} if
$$\ell([g,h]) \leq 4 \cdot \ell(g) \ell(h), \quad \forall g,h \in G.$$
\end{definition}

Commutator-contractive length functions naturally arise on unitary groups of $C^*$-algebras. Indeed, the following lemma shows that for a unital $C^*$-algebra $A$, the function $\ell(g) := \frac12\|1 - g\| \in [0,1]$ is a commutator-contractive invariant length function on $U(A)$.

\begin{lemma} \label{comp} Let $A$ be unital $C^*$-algebra.
Let $g,h \in U(A)$. The following relations hold for all $g,h \in U(A)$:
$$\ell(gh) \leq \ell(g) + \ell(h), \quad \ell(g) = \ell(\bar g), \quad \ell(gh)=\ell(hg) \quad \mbox{and} \quad \ell([g,h]) \leq 4 \cdot \ell(g) \ell(h).$$ 
\end{lemma}
\begin{proof}
The last inequality is the only non-trivial assertion. Let us compute:
$$\ell([g,h])=\frac{\|1 - gh\bar g \bar h\|}2 = \frac{\|hg -gh\|}2 = \frac{\|(1-h)(1-g) - (1-g)(1-h)\|}2 \leq 4 \cdot \ell(g) \ell(h).$$
This finishes the proof.
\end{proof}

\begin{remark}
Note that the usual normalized Hamming metric on the symmetric group $S_n$, which is given by the formula
$$d(\rho,\sigma) = \frac{|\{i \mid \rho(i) \neq \sigma(i)\}|}n,$$
is not commutator-contractive. However, we have an analogous inequality
$$\ell([\rho,\sigma]) \leq 2 \cdot \min\{\ell(\rho),\ell(\sigma)\}, \quad \forall \sigma,\rho \in S_n$$
which holds in every group with an invariant length function.
\end{remark}

We will now explain what we mean by approximation. We will state a definition and two propositions which clarify various aspects.

\begin{definition}
Let $\cC$ be a class of groups with an invariant length function.
A group $G$ is said to have the $\cC$-\emph{approximation property} if the following holds. For all $g \in G \setminus \{e\}$, there exists $\delta_g >0$ such that
for all finite subsets $F \subset G$, and $\varepsilon>0$, there exists a group $C \in \cC$ and a map $\varphi \colon G \to C$ such that
\begin{enumerate}
\item $\varphi(e)=e$,
\item $d(\varphi(gh),\varphi(g)\varphi(h))< \varepsilon$, for all $g,h \in F$,
\item $\ell(\varphi(g)) \geq \delta_g$, for all $g \in F \setminus \{e\}$.
\end{enumerate}
\end{definition}

The following results give easy reformulations of the above definition. Proofs of comparable results can be found in \cite{MR2455513} and \cite{MR2178069}.

\begin{proposition} \label{propref1} Let $G$ be a countable group and let $\pi \colon F \to G$ be a surjection from a countable free group. Let $\cC$ be a class of groups with an invariant length function. The group $G$ has the $\cC$-approximation property if and only if there exists a sequence of groups $C_n \in \cC$ and a sequence of homomorphisms $\varphi_n \colon F \to C_n$ such that
$$ \lim_{n \to \infty} \ell(\varphi_n(g)) = 0, \quad \forall g \in \ker(\pi),$$
and
$$ \liminf_{n \to \infty} \ell(\varphi_n(g)) \neq 0, \quad \forall g \not\in \ker(\pi).$$
\end{proposition}

\begin{proposition}
Let $G$ be a countable group and let $\cC$ be a class of groups with an invariant length function. The group $G$ has the $\cC$-approximation property if and only if there exists a non-principal ultrafilter $\omega$ on $\N$, a sequence $C_n$ of groups in $\cC$, and an injective homomorphism
$$\iota \colon G \hookrightarrow \frac{\prod_{n \in \N} C_n}{ \left\{ (g_n)_{n \in \N} \mid \ell_{C_n}(g_n) \to_{\omega} 0 \right\} }.$$ 
\end{proposition}

Note that the group $\frac{\prod_{n \in \N} C_n}{ \left\{ (g_n)_{n \in \N} \mid \ell_{C_n}(g_n) \to_{\omega} 0 \right\} }$ can be seen as a metric ultraproduct of groups with an invariant length function. Note that everything fits well with model theory for metric structures, an extension of the usual model theory which was developed in \cite{MR2436146}.

\vspace{0.2cm}

In the literature, various classes of groups have been considered. The embedability of the group von Neumann algebra of a group $G$ into an ultraproduct of hyperfinite II$_1$-factors is equivalent (see \cite{MR2436761}) to the $\cU_{\rm HS}$-approximation property of $G$, where $\cU_{\rm HS}$ is the class of finite-dimensional unitary groups with the normalized Hilbert-Schmidt metric. Following \cite{MR2436761}, such groups are called \emph{hyperlinear}. The MF-embedability of $C^*(G)$ or $C^*_{\rm red}(G)$ in the sense of Blackadar-Kirchberg (see \cites{MR1437044, MR1796908, MR1835511}) implies the $\cU_{\|.\|}$-approximation property of $G$, where $\cU_{\|.\|}$ is the class of finite-dimensional unitary groups with the normalized operator norm. It is conjectured that wide classes of finite $C^*$-algebras -- in particular, all reduced group $C^*$-algebras -- do satisfy MF-embedability. Note that the length functions in the class $\cU_{\|.\|}$ are commutator-contractive.

A group $G$ is \emph{sofic} if and only if $G$ has the $\cS$-approximation property, where $\cS$ denotes the class of symmetric groups with the usual normalized Hamming distance. A group is said to be \emph{weakly sofic} iff it has the $\cF$-approximation property, where $\cF$ is the class of all finite groups with an invariant length function, see \cite{MR2455513}. We also consider the class $\cF_{\rm c}$ of finite groups with a commutator-contractive invariant length function. Note that finite subgroups of unitary groups with the normalized operator norm belong to $\cF_{\rm c}$, so that there is an abundance of finite groups with commutator-contractive invariant length functions. A group is LEF (locally embedable into finite groups, see \cite{MR1458419}) if it has the $\cF_{\delta}$-approximation property, where $\cF_{\delta}$ denote the class of finite groups with the discrete metric, i.e. $\ell(g)=1$, for all $g \neq e$. This last example is trivial from the metric point of view and it is well-known that there are many groups which are not LEF, see \cite{MR1458419}. In \cite{MR2566306}, the author constructed a hyperlinear group with Kazhdan's property (T), which is not LEF.

\vspace{0.2cm}

For the time being, there is nothing known about the existence of groups which do \emph{not} have the approximation properties with respect to the classes $\cU_{\rm HS}$, $\cU_{\|.\|}$, $\cS$ or $\cF$. Hence, it seems acceptable and fair to try to put even more conditions on the class of groups in order to be able to obtain at least some positive results.

The most restrictive non-trivial class of groups with an invariant metric seems to be the class $\cF_{\rm c}$ of finite groups with a commutator-contractive invariant length function.
In this note, we want to develop the basic theory of groups in the class $\cF_{\rm c}$ and prove that Graham Higman's group \cite{MR0038348} does not have the $\cF_{\rm c}$-approximation property. Higman's group was the first example of a finitely presented group without any finite quotients. Hence, it seems to be a natural candidate to consider when searching for a non-sofic or non-hyperlinear group.

To the best knowledge of the author, the class of groups with the $\cF_c$-approximation property is not known to be bigger than the class of LEF groups; in particular, it is \emph{not} known whether all amenable groups have the $\cF_c$-approximation property. It would also be interesting do decide whether the $\cU_{\|.\|}$-approximation property implies the $\cF_c$-approximation property.

\section{Finite groups with an invariant length function}

Throughout this section, let $G$ be a finite group with an invariant length function $\ell_G \colon G \to [0,1]$. If $H \subset G$ is a normal subgroup, then $G/H$ carries an induced invariant length function
$$\ell_{G/H}(gH) = \min \{ \ell_G(k) \mid k \in gH \}.$$
Note that the natural projection $\pi \colon G \to G/H$ is a contraction.
If $\ell_G$ was commutator-contractive, then $\ell_{G/H}$ is commutator-contractive as well. Indeed, let $g,k \in G$ be representatives of cosets $gH$ and $kH$ of minimal length, then
$$\ell_{G/H}([g,k]H) = \min \{ \ell_G(g') \mid g' \in [g,k]H \} \leq \ell_G([g,k]) \leq 4 \cdot \ell_{G/H}(gH) \ell_{G/H}(kH).$$

\begin{lemma} \label{comp}
Let $g,h,k \in G$ and let $\ell\colon G \to [0,1]$ be a commutator-contractive invariant length function on $G$. Then,
$$d([g,h],[g,k]) \leq 4 \cdot d(h,k) \ell(g).$$
\end{lemma}
\begin{proof}
We compute
$$d([g,h],[g,k]) = \ell(gh\bar g \bar h k g \bar k \bar g) 
= \ell(\bar kh\bar g \bar h k g ) 
= \ell([\bar kh, \bar g]) 
 \leq  4 \cdot \ell(\bar k h) \ell(g) 
= 4 \cdot d(h,k) \ell(g).
$$
This finishes the proof.
\end{proof}

We set $\delta(G) := \min \{ \ell(g) \mid g \in G, g \neq e \}.$ A group with an invariant length function is said to be \emph{discrete} if $\delta(G)=1$.
We define $G_{\varepsilon}$ to be the subgroup which is generated by the set $\{g \in G \mid \ell(g) \leq \varepsilon \}$ and set
$\eta(G) := \min \{ \varepsilon \in [0,1] \mid G_{\varepsilon} = G \}.$
Note that $G_{\varepsilon} \subset G$ is automatically a normal subgroup.

\begin{lemma} \label{lem1}
Let $G$ be a finite group with an invariant length function.
Then, either
$\delta(G/G_{\varepsilon}) > \varepsilon$ or $G_{\varepsilon}=G$.
\end{lemma}
\begin{proof}
Let $g$ be such that $\ell_{G/G_{\varepsilon}}(gG_{\varepsilon}) \leq \varepsilon$. Then there exists $h \in G_{\varepsilon}$ such that $\ell_G(gh) \leq \varepsilon$. Hence $g \in G_{\varepsilon}$ and $gG_{\varepsilon}=G_{\varepsilon}$. This implies $\delta(G/G_{\varepsilon}) > \varepsilon$
unless $G_{\varepsilon}=G$.
\end{proof}
\begin{lemma} \label{lem2}
Let $G$ be a finite group with an invariant length function.
If $\eta(G)> \varepsilon$, Then
$$\eta(G/G_{\varepsilon}) = \eta(G).$$
\end{lemma}
\begin{proof}
The inequality $\eta(G/G_{\varepsilon}) \leq \eta(G)$ is obvious since the projection $\pi \colon G \to G/G_{\varepsilon}$ is a contraction. If $G/G_{\varepsilon}$ is generated by elements of length $\eta(G/G_{\varepsilon})$, then $G_{\eta(G/G_{\varepsilon})}$ together with $G_{\varepsilon}$ will generate $G$. Hence, we get:
$$\eta(G) \leq \max\{\varepsilon, \eta(G/G_{\varepsilon}) \} \leq \eta(G).$$
Here, the last inequality follows from our assumption $\varepsilon< \eta(G)$ and the observation $\eta(G/G_{\varepsilon}) \leq \eta(G)$. Since $\varepsilon< \eta(G)$, we must have $\eta(G/G_{\varepsilon}) = \eta(G)$. This finishes the proof.
\end{proof}

\begin{remark} \label{subgr}
We also note that $\eta(G_{\varepsilon}) \leq \varepsilon$. Similarly, $\delta(G_{\varepsilon}) = \delta(G)$ if $\delta(G) \leq \varepsilon$.
\end{remark}

\begin{definition} Let $G$ be a group. The lower central series $(\gamma_n(G))_{n \in \N}$ is defined as $\gamma_0(G) = G$ and $\gamma_{n+1}(G) = [\gamma_n(G),G]$. A group $G$ is said to be nilpotent if $\gamma_n(G)= 1$ for some $n \in \N$. The nilpotency class of $G$ is defined to be ${\rm nil}(G) := \min\{n \in \N \mid \gamma_n(G) = 1 \}$.
\end{definition}

Note that $G_{\delta(G)}$ is abelian if $\delta(G)<1/4$. Indeed, $\ell([g,h]) \leq 4 \ell(g)\ell(h) < \delta(G)$ if $\ell(g)=\ell(h)=\delta(G)$. More generally, $G_{\varepsilon}$ is nilpotent if $\varepsilon$ is sufficiently small.
Zassenhaus' Lemma states that a discrete subgroup of a Lie group which is generated by elements which are sufficiently close to $1$ is nilpotent. The following Proposition is the quantitative analogue in our situation. This result is not needed in the proof of our main results, but of independent interest.

\begin{proposition}[Zassenhaus' lemma] \label{zas}
Let $G$ be a finite group with a commutator-contractive invariant length function $\ell \colon G \to [0,1]$. Let $\varepsilon< 1/4$.  
The group $G_{\varepsilon}$ is nilpotent with
$$ {\rm nil}(G_{\varepsilon}) \leq \frac{\ln(4\delta(G))}{\ln(4 \varepsilon)}.$$  
\end{proposition}
\begin{proof}
If $\delta(G) \geq 1/4$, then $G_{\varepsilon}=1$ and there is nothing to prove.
Since $\ell([g,h]) \leq 4 \ell(g) \ell(h)$ for all $g,h \in G$, we conclude from standard commutator identities that
$[G_{\varepsilon'},G_{\varepsilon}] \subset G_{4 \varepsilon \varepsilon'}$ for all $\varepsilon'>0$.
Hence, we see by induction that $\gamma_n(G_{\varepsilon}) \subset G_{(4\varepsilon)^n \varepsilon}$. Since $G_{\delta'}=1$ for $\delta' < \delta(G)$, we conclude that
$\gamma_n(G_{\varepsilon}) = 1$ if $(4\varepsilon)^n \varepsilon < \delta(G)$. Equivalently, $\gamma_n(G_{\varepsilon}) = 1$ if
$$n > \frac{\ln(\delta(G)) - \ln(\varepsilon)}{\ln(4 \varepsilon)}.$$ 
In particular, we get
$$ {\rm nil}(G_{\varepsilon}) \leq 1 + \frac{\ln(\delta(G)) - \ln(\varepsilon)}{\ln(4 \varepsilon)} = \frac{\ln(4\delta(G))}{\ln(4 \varepsilon)}.$$  
This finishes the proof.
\end{proof}

\begin{corollary} \label{corzas}
Let $G$ be a finite group with a commutator-contractive invariant length function $\ell \colon G \to [0,1]$. If $\eta(G) < 1/4$, then
$${\rm nil}(G) \leq \frac{\ln (4 \delta(G))}{\ln(4 \eta(G))}.$$
\end{corollary}
\begin{proof}
By definition of $\eta(G)$, we have $G_{\eta(G)} = G$. The claim follows from Proposition \ref{zas}.
\end{proof}

%\section{Examples}
%
%Let $(X,d)$ be a metric space with ${\rm diam}(X) = 1$ and let $G \times X \to X$, $(g,x) \mapsto g.x$ be an isometric $G$-action. We set
%$$\ell(g) = \sup_{x \in X} d(g.x,x).$$
%
%
%
%\section{Only questions}
%
%Let $\phi \colon \F_n \to U(n)$ be a homomorphism and set $G:= \varphi(\F_n)$ with the restricted length function. Consider the subgroups $G_{\varepsilon} \subset G$.
%
%\begin{enumerate}
%\item $G_{\varepsilon}$ is a closed and open subgroup of $G$. Indeed, if $d(g,G_{\varepsilon}) \leq \varepsilon$, then $g \in G_{\varepsilon}$.
%\item $G_{\varepsilon}$ has finite index in $G$.
%\item We set
%$$\underline{\ell}(g) := \lim_{\varepsilon \to 0} \ell_{G/G_{\varepsilon}} (gG_{\varepsilon}) \leq \ell_G(g).$$
%\end{enumerate}
%
%

\section{The main result}

Higman's group as defined in \cite{MR0038348} is given by a concrete presentation:
$$H = \langle a_0,a_1,a_2,a_3 \mid a_i = [a_{i+1},a_i], \forall i \in \Z/4\Z \rangle.$$
It was the first example of a finitely presented group without any finite quotients. We want to show that homomorphisms from the free group on four generators to a finite group with a commutator-contractive invariant length function, such the the images of the generators of the free group satisfy the relations of Higman's group up to some $\varepsilon$ are either trivial or uniformly non-trivial on generators. Trivial means in this context, that the generators of the free group are mapped to elements whose length is comparable with $\varepsilon$. Later, we will exclude the possibility of the second case. Our first result is:

\begin{proposition} \label{mainprop}
Let $G$ be a finite group with a commutator-contractive invariant length function $\ell \colon G \to [0,1]$. Let $\varepsilon < 1/64$.
Let $\varphi \colon \F_4 = \langle b_0,b_1,b_2,b_3 \rangle \to G$ be a surjective homomorphism and set $a_i := \varphi(b_i)$. Assume that
$$d(a_i,[a_{i+1},a_i]) \leq \varepsilon, \quad \forall i \in \Z/4\Z.$$
Then, either $\ell(a_i) < 4\varepsilon$, for all $i \in \Z/4\Z$, or
$\ell(a_i) \geq 7/32$, for all $i \in \Z/4\Z$.
\end{proposition}
\begin{proof}
We define $${\rm max}(\varphi) := \max\{ \ell(a_i) \mid i \in \Z/4\Z \} \quad \mbox{and} \quad
{\rm min}(\varphi) := \min\{ \ell(a_i) \mid i \in \Z/4\Z \}.$$ We clearly have ${\rm min}(\varphi) \leq {\rm max}(\varphi)$.
We have to show that, if there exists $a_i$ such that $\ell(a_i) \geq 4\varepsilon$, then $\ell(a_i) \geq 7/32$ for all $i \in \Z/4\Z$.
Since
$$\ell(a_i) - \varepsilon \leq \ell([a_{i+1},a_i]) \leq 4 \cdot \ell(a_{i+1}) \ell(a_i)$$
we conclude
\begin{equation} \label{eqest}
\frac3{16} \leq \frac14 - \frac{\varepsilon}{ 4 \ell(a_i)} \leq \ell(a_{i+1}).
\end{equation}
Hence, $\ell(a_{i+1}) \geq \frac3{16} \geq 4\varepsilon$ and after completing the cycle, we obtain $\ell(a_i) \geq 3/16$ for all $i \in \Z/4$. Using Equation (\ref{eqest}) again, we get with $\varepsilon < 1/64 \leq 3/128$ that
$$\frac{7}{32} = \frac14 - \frac{1}{32} \leq \frac14 - \frac{4 \varepsilon}{3}= \frac14 - \frac{\varepsilon}{4 \cdot 3/16} \leq \frac14 - \frac{\varepsilon}{ 4 \ell(a_{i-1})} \leq \ell(a_i)$$ for all $i \in \Z/4 \Z$.
\end{proof}

\begin{theorem} \label{main}
Let $G$ be a finite group with a commutator-contractive invariant length function $\ell \colon G \to [0,1]$. Let $\varepsilon < 1/64$.
Let $\varphi \colon \F_4 = \langle b_0,b_1,b_2,b_3 \rangle \to G$ be a surjective homomorphism and set $a_i := \varphi(b_i)$. Assume that
$$d(a_i,[a_{i+1},a_i]) \leq \varepsilon, \quad \forall i \in \Z/4\Z.$$
Then, $\ell(a_i) < 4 \varepsilon$ for all $i \in \Z/4\Z$.\end{theorem}
\begin{proof}
First of all, we may assume that $\delta(G)<1/64$, since otherwise $\varphi$ is well-defined on $H$, but there are no finite quotients of $H$.

Let us fix some $\varepsilon<1/64$ and argue with proof by contradiction.
Let  $\varphi \colon \F_4 \to G'$ be a homomorphism such that
$d(a_i,[a_{i+1},a_i]) \leq \varepsilon$, for all $i \in \Z/4 \Z$ and $\ell(a_i) \geq 4 \varepsilon$ for some $i \in \Z/4\Z$. We set $n:= |G'|$ and consider the set of all homomorphisms $\varphi\colon \F_4 \to G$ satisfying the above conditions with $|G| \leq n$. By compactness, there exists $\varphi \colon \F_4 \to G$ as above with $\delta(G)$ maximal.

Let us first assume that $G_{\delta(G)} \neq G$.
Since then $\delta(G/G_{\delta(G)}) > \delta(G)$ by Lemma \ref{lem1}, the induced homomorphism $\bar \varphi \colon \F_4 \to G/G_{\delta(G)}$ violates one of the conditions. Since the projection $\pi \colon G \to G/G_{\delta(G)}$ is contractive, the only condition that can be violated is the lower bound on $\ell(a_i)$. Hence, there exists $i \in \Z/4\Z$ such that $\ell_{G/G_{\delta(G)}}(a_i) < 4 \varepsilon < 7/32$.

Using Proposition \ref{mainprop}, we can conclude that
$\ell_{G/G_{\delta(G)}}(a_i) < 4 \varepsilon$ for all $i \in \Z/4\Z$. 
Hence, we find
$\tilde{a}_i \in G_{\delta(G)}$, for $i \in \Z/4\Z$, such that $d(\tilde a_i,a_i) < 4 \varepsilon$. Obviously, the same is true if $G_{\delta(G)}=G$. Since $\delta(G) < 1/4$, we know that $G_{\delta(G)}$ is abelian. We compute, using Lemma \ref{comp}, that for all $i \in \Z/4\Z$:
\begin{eqnarray*}
\ell(a_i) &\leq& \ell([a_{i+1},a_i]) + \varepsilon \\
& \leq & 4 \varepsilon \cdot \ell(a_{i+1}) + \ell([a_{i+1},\tilde a_i]) + \varepsilon \\
& \leq & 4 \varepsilon \cdot (\ell(a_{i+1})+ \ell(\tilde a_{i})) + \ell([\tilde a_{i+1},\tilde a_i]) + \varepsilon  \leq 9 \varepsilon.\\  
\end{eqnarray*}
Since $9\varepsilon < 7/32$, we conclude from Proposition \ref{mainprop} that $\ell(a_i) < 4 \varepsilon$ for all $i \in \Z/4\Z$. This is a contradiction to our assumption and finishes the proof.
\end{proof}

\begin{corollary} \label{corhig}
Higman's group does not have the $\cF_{\rm c}$-approximation property. 
\end{corollary}
\begin{proof}
We are using the reformulation of the $\cF_{\rm c}$-approximation property according to Proposition \ref{propref1}.
Suppose that $H$ has the $\cF_{\rm c}$-approximation property with certain $\delta_{a_i}>0$, for $i \in \Z/4\Z$, and let $0<\varepsilon< 1/3000$ be arbitrary. Then, there exists a finite group $G$ with a commutator-contractive invariant length function, and a homomorphism $\varphi \colon \F_4 \to G$ such that
$d(a_i,[a_{i+1},a_i]) \leq \varepsilon$, for all $i \in \Z/4 \Z$ and $\ell(a_i) \geq \delta_{a_i}/2$ for all $i \in \Z/4\Z$. By the preceding theorem,
$\delta_{a_i}/2 < 56 \varepsilon$ which is impossible for $\varepsilon$ small enough. This finishes the proof.
\end{proof}

\begin{remark}
Theorem \ref{main} implies with a similar argument as in Corollary \ref{corhig} that any quotient of Higman's group does not have the $\cF_{\rm c}$-approximation property.
\end{remark}

%\section{Some remarks}
%
%Let $G$ be a finitely presented group and $\phi \colon \F_n \to G$ be a surjection. We denote by $a_1,\dots,a_n \in G$ the images of the generators of $\F_n$. Assume that $G$ is normally generated by $g \in G$. Let $r_1,\dots,r_k \in \F_n$ be a finite set of defining relations of $G$, i.e.\ 
%$$G = \langle g_1,\dots,g_n \mid r_1,\dots,r_k \rangle.$$
%
%Let $n$ be some integer, such that each of the $a_i$'s is a product of $n$ conjugates of $r_i$'s and the element $g$. Let us fix $0 < \varepsilon < 1/4$. Consider a homomorphism $\psi \colon \F_n \to (H,\ell_{H})$ such that
%$\ell_H(\psi(r_i)) \leq \varepsilon$ for all $1 \leq i \leq k$, and assume that $\ell_H(\psi(g)) \geq \varepsilon$. By compactness, we may assume that $\psi \colon \F_n\to (H,\ell_H)$ is such that $\delta(H)$ is maximal among all such $\psi$ with a fixed $\varepsilon$ and groups of size at most $|H|$.
%
%\begin{enumerate}
%\item If $\delta(H) > \varepsilon$, then $\psi$ defines a homomorphism from $G$ to $H$.
%\item If $\delta(H)< \varepsilon < 1/4$, then $H_{\delta(H)}$ is abelian and $\delta(H/H_{\delta(H)}) > \delta(H)$. In particular, $\ell_{H/H_{\delta(H)}}(g H_{\delta(H)})< \varepsilon$.
%\end{enumerate}
%

\end{document}